\documentclass[11pt]{article}

\usepackage{IEEEtrantools}
\usepackage{latexsym}
\usepackage{amsfonts}
\usepackage{amssymb}
\usepackage{amsmath,amsfonts,amsthm}
\usepackage{mathrsfs}
\usepackage{multirow}
\usepackage{tikz}
 \usetikzlibrary{arrows,decorations.markings}

\newcommand{\be}{\begin{equation}}
\newcommand{\ee}{\end{equation}}
\newcommand{\bea}{\begin{eqnarray}}
\newcommand{\eea}{\end{eqnarray}}
\newcommand{\bean}{\begin{eqnarray*}}
\newcommand{\eean}{\end{eqnarray*}}
\newcommand{\brray}{\begin{array}}
\newcommand{\erray}{\end{array}}
\newcommand{\biearray}{\begin{IEEEarray}{rCl}}
\newcommand{\eiearray}{\end{IEEEarray}}

\newcommand{\newsection}[1]{\setcounter{equation}{0}
\setcounter{dfn}{0}
\section{#1}}

\newtheorem{dfn}{Definition}[section]
\newtheorem{thm}[dfn]{Theorem}
\newtheorem{lmma}[dfn]{Lemma}
\newtheorem{ppsn}[dfn]{Proposition}
\newtheorem{crlre}[dfn]{Corollary}
\newtheorem{xmpl}[dfn]{Example}
\newtheorem{rmrk}[dfn]{Remark}

\newcommand{\bdfn}{\begin{dfn}\rm}
\newcommand{\bthm}{\begin{thm}}
\newcommand{\blmma}{\begin{lmma}}
\newcommand{\bppsn}{\begin{ppsn}}
\newcommand{\bcrlre}{\begin{crlre}}
\newcommand{\bxmpl}{\begin{xmpl}}
\newcommand{\brmrk}{\begin{rmrk}\rm}

\newcommand{\edfn}{\end{dfn}}
\newcommand{\ethm}{\end{thm}}
\newcommand{\elmma}{\end{lmma}}
\newcommand{\eppsn}{\end{ppsn}}
\newcommand{\ecrlre}{\end{crlre}}
\newcommand{\exmpl}{\end{xmpl}}
\newcommand{\ermrk}{\end{rmrk}}

\newcommand{\bbc}{\mathbb{C}}

\newcommand{\bbn}{\mathbb{N}}
\newcommand{\bbr}{\mathbb{R}}

\newcommand{\cla}{\mathcal{A}}

\newcommand{\clh}{\mathcal{H}}

\newcommand{\cll}{\mathcal{L}}

\newcommand{\cls}{\mathcal{S}}

\newcommand{\prf}{\noindent{\it Proof\/}: }

\def \qed { \mbox{}\hfill
$\Box$\vspace{1ex}}

\begin{document}
 
 
 
\author{\sc  { Bipul Saurabh}}
\title{ Spectral dimension of quaternion spheres}
\maketitle
 \begin{abstract}
Employing ideas of noncommutative geometry, certain dimensional invariant for quantum homogeneous spaces has been 
proposed and here we take up its computation for   quaternion spheres.
\end{abstract}
{\bf AMS Subject Classification No.:} {\large 58}B{\large 34}, {\large 58}B{\large 32}, {\large46}L{\large 87} 

\newsection{Introduction}  
In Connes' (\cite{Con-1994aa}) formulation, a geometric space is given by a triple $(\cla,\clh,D)$ called a spectral triple, with $\cla$ being an involutive algebra represented as bounded operators
on a Hilbert space $\clh$, and $D$
being a selfadjoint operator with compact resolvent and having
bounded commutators with the algebra elements. Connes further defined the dimension of a spectral triple to be the quantity $\inf\{\delta: \mbox{Tr}(|D|^{-\delta})<\infty \}$. Utilizing this,  Chakraborty and Pal
(\cite{ChaPal-2016aa}) introduced an invariant called spectral dimension, for an ergodic $C^*$-dynamical system or equivalently for a homogeneous space of a compact quantum group. They considered 
all finitely summable equivariant  spectral triples on the GNS space of the  state invariant under the group action and defined the spectral dimension of the homogeneous space  to be the infimum of
the summability of the associated Dirac operators. 
Compact quantum groups and their quotient spaces   are natural examples of homogeneous spaces.  
In the same paper \cite{ChaPal-2016aa}, Chakraborty and Pal computed spectral dimension of many such homogeneous spaces, both in classical
and quantum situations and   it was conjectured that the spectral dimension of a homogeneous space of a (classical) compact Lie group
is same as its dimension as a differentiable manifold. The spectral dimensions of  $SU(2)$, $SU_q(n)$ and $S_q^{2n+1}$ point towards this conjecture. All these 
examples are homogeneous spaces of type $A$ quantum groups. 
To examine the conjecture, we need to explore more examples.  In this article, we  take up 
the case of 
quaternion sphere $H^n$. These spaces are homogeneous spaces of type $C$ quantum groups.  
Here we show that  the spectral dimension is  equal to its dimension as a real manifold. Therefore it strengthens 
the conjecture of  Chakraborty and Pal. The computation of the invariant for $H^n$ is the first instance of computation 
for  homogeneous spaces of type $C$ quantum groups.

We will sometimes write a spectral triple $(\cla,\clh,D)$ 
as $(\clh,\pi,D)$ where $\pi$ is the representation of $\cla$ in the Hilbert space $\clh$. For a subset $S$ of a $C^*$-algebra, 
$\overline{S}$ will denote the closed linear span of $S$ in $A$.

\newsection{Spectral dimension}
In this section, we recall from \cite{ChaPal-2016aa} the definition of spectral dimension of a $C^*$-dynamical system.
Let us begin with the definition of a homogeneous space. 
\bdfn A compact quantum group $G$ acts on a $C^*$-algebra $A$ if there exists a $*$-homomorphism $\tau: A \rightarrow A \otimes C(G)$ such that 
\begin{enumerate}
 \item $(\tau \otimes id)\tau=(id\otimes \Delta)\tau$,
 \item$\overline{\{(I\otimes b)\tau(a): a\in A, b \in C(G)\}}=A\otimes C(G)$.
\end{enumerate}
where $\Delta$ is the comultiplication map of $G$.
We call an action $\tau$ homogeneous or ergodic if the fixed point subalgebra $\{a \in A : \tau(a)=a\otimes I\}$ is $\bbc I$. In that case, the associated $C^*$-algebra $A$ is called an homogeneous space 
of $G$ and the triple $(A,G,\tau)$ is called an ergodic $C^*$-dynamical system.
\edfn
A covariant representation of a $C^*$-dynamical system $(A,G,\tau)$ is a pair $(\pi,U)$ consisting of a representation $\pi:A \rightarrow \cll(\clh)$ and a unitary representation 
of $G$ on $\clh$ such that for all $a \in A$, one has 
\[
 (\pi \otimes id)\tau(a)= U(\pi(a) \otimes I)U^*.
\]
\bdfn
Let $(\pi,U)$ be a covariant representation of a $C^*$-dynamical system $(A,G,\tau)$ and $(\clh,\pi,D)$ be a spectral triple 
for a dense $*$-subalgebra $\cla$ of $A$. We call $(\clh,\pi,D)$ equivariant with respect to $(\pi,U)$ if $D\otimes I$ commutes with $U$.
\edfn

Associated with a homogeneous action $\tau$ of $G$ is a unique invariant state $\rho$ on the homogeneous space $A$ that obeys 
\[
 (\rho \otimes id)\tau(a)= \rho(a)I, \qquad a\in A.
\]
Consider the GNS representation $(\clh_{\rho}, \pi_{\rho}, \eta_{\rho})$ of $A$ associated with the state $\rho$. Using the invariance property of $\tau$, one can show that the 
action $\tau$ induces a unitary representation $U_{\tau}$ of $G$ on $\clh_{\rho}$ and the pair $(\pi_{\rho}, U_{\tau})$ is a covariant representation of the
system $(A,G,\tau)$. Let $\mathcal{O}(G)$ be the dense $*$-Hopf subalgebra of $C(G)$ generated by matrix entries of  irreducible unitary representations 
of $G$. Define 
\[
 \cla:= \{ a \in A: \tau(a) \in A \otimes_{alg}\mathcal{O}(G)\}.
\]
It follows from part $(1)$ of Theorem $1.5$ in \cite{Pod-1995aa} that $\cla$ is a dense $*$-subalgebra of $A$. Define $\xi$ to be the class of spectral triples of $\cla$ equivariant
 with respect to the covariant representation $(\pi_{\rho}, U_{\tau})$. The spectral dimension denoted by $\cls dim(A,G,\tau)$ of the $C^*$-dynamical system $(A,G,\tau)$ is defined to be the quantity
 \[
  \inf\{p>0: \exists  D \mbox{ such that }  (\cla,\clh_{\rho},D) \in \xi \mbox{ and } D \mbox{ is $p$-summable}\}.
 \]

\newsection{Main result}
Here we briefly recall some notions related to  quaternion spheres  $H^n$ (or $SP(2n)/SP(2n-2)$)  and then compute its spectral dimension. For $1 \leq i,j \leq 2n$, define a continuous map 
\begin{IEEEeqnarray}{rCl}
 u_i^j:SP(2n) & \rightarrow & \bbc;  \quad A \mapsto a_j^i \nonumber
\end{IEEEeqnarray} 
where $a_j^i$ is the $ij^{th}$ entry of $A \in SP(2n)$. 
The $C^*$-algebra  $C(SP(2n))$ is generated by elements of the set $\{u_{j}^{i}:1\leq i,j\leq 2n\}$. 
In the same way, define the generators  $\{v_{j}^{i}:1\leq i,j\leq 2n-2\}$  of $C(SP(2n-2))$. 
Define the map $\Phi: C(SP(2n)) \rightarrow C(SP(2n-2))$ 
as follows.
\[
\Phi(u_{j}^{i})=\begin{cases}
                             v_{j-1}^{i-1}, 
                             & \mbox{ if } i \neq 1 \mbox{ or } 2n, \mbox{ or } j \neq 1 \mbox{ or } 2n, \cr
			 \delta_{ij}, & \mbox{ otherwise. } \cr
			\end{cases}
\] 
Clearly $\phi$ is a $C^{*}$-epimorphism  obeying $\Delta\phi=(\phi\otimes \phi)\Delta$ where $\Delta$ is the co-multiplication map of $C(SP(2n))$.
In such a case, one defines the quotient space $C(SP(2n)/SP(2n-2))$ by,
\[
         C(SP(2n)/SP(2n-2)) = \left\{a\in C(SP(2n)) : (\phi\otimes id)\Delta(a) = I\otimes a\right\}. 
\] 
The quotient space $SP(2n)/SP(2n-2)$ can be realized as the $n$-dimensional quaternion sphere $H^{n}$.  Also, each of 
the generators $\{u_{j}^{1}:   1\leq j\leq 2n\}$ can be viewed as 
projection on to a fixed complex coordinate of a point in $H^{n}\subset \bbc^{2n}$ and for  $ 1\leq j\leq 2n$, the map $u_{j}^{2n}$ is the complex conjugate of $u_{2n+1-j}^{1}$. 
This shows that  $C(SP(2n)/SP(2n-2))$ is generated by $\{u_{j}^{i}: i=1 \mbox{ or } 2n, 1\leq j\leq 2n\}$.
Restricting the co-multiplication map to $C(SP(2n)/SP(2n-2))$  gives an action $\tau$ of  the compact quantum group $SP(2n)$ on $SP(2n)/SP(2n-2)$. 
\begin{IEEEeqnarray}{rCl}
 \tau: C(SP(2n)/SP(2n-2)) &\longrightarrow & C(SP(2n)/SP(2n-2))\otimes C(SP(2n)) \nonumber \\
 a &\longmapsto & \Delta a \nonumber .
\end{IEEEeqnarray}
It is not difficult to verify that the system $(C(SP(2n)/SP(2n-2)), SP(2n), \tau)$ is an ergodic $C^*$-dynamical system and the invariant state $\rho$ of $\tau$ is the  faithful Haar state $h$ of 
$C(SP(2n))$ restricted to $C(SP(2n)/SP(2n-2))$. 
By Theorem $1.5$ of \cite{Pod-1995aa}, we get
\begin{IEEEeqnarray}{rCl}\label{eq1}
C(SP(2n)/SP(2n-2))= \overline{\oplus_{\lambda \in \widehat{SP(2n)}}\oplus_{i \in I_{\lambda} }W_{\lambda,i}}
\end{IEEEeqnarray}
where $\lambda=(\lambda_1,\lambda_2,\cdots , \lambda_n)$ represents
the highest weight of  a finite-dimensional irreducible co-representation $u_{\lambda}$ of $C(SP(2n))$,  $I_{\lambda}$ is
the multiplicity of $u_{\lambda}$ 
and $ W_{\lambda, i}$ corresponds to 
$u_{\lambda}$ in the sense of Podles (see page $4$, \cite{Pod-1995aa}) for all $i \in I_{\lambda}$. Using Zhelobenko branching rule (see page 79, \cite{Zhe-1962aa} 
and Theorem 1.7 in \cite{Pod-1995aa}, page 145-146 in \cite{Sau-2015aa}), we get  
\begin{displaymath}
I_{\lambda}= \begin{cases}
                             \lambda_1-\lambda_2+1, & \quad \mbox{if} \quad \lambda_i=0 \quad \mbox{for all}\quad  i \geq 3, \cr
                             0, & \quad \mbox{otherwise}. \cr
                            \end{cases}
\end{displaymath}
Define 
\[
 \mathcal{O}(SP(2n)/SP(2n-2)):=\oplus_{\lambda \in \widehat{SP(2n)}}\oplus_{i \in I_{\lambda} }W_{\lambda,i}. 
\]
It is not difficult to see that 
$\mathcal{O}(SP(2n)/SP(2n-2))$ is the algebra generated by the elements of the set $\{u_{j}^{i}: i=1 \mbox{ or } 2n, 1\leq j\leq 2n\}$.
Moreover,  the algebra $\mathcal{O}(SP(2n)/SP(2n-2))$ is a dense Hopf $*$-algebra consisting of all $a \in C(SP(2n)/SP(2n-2))$ such that $\tau(a)\in C(SP(2n)/SP(2n-2))\otimes_{alg}\mathcal{O}(SP(2n))$.

Let  $U(\mathfrak{sp}(2n))$ be the universal enveloping algebra of the  Lie algebra $\mathfrak{sp}(2n)$. 
We will view $\mathfrak{sp}(2n)$ as a subset of $U(\mathfrak{sp}(2n))$. Then $U(\mathfrak{sp}(2n))$ is 
generated by  $H_i,E_i, F_i \in  \mathfrak{sp}(2n)$,  $i=1,2,\cdots,n$, 
satisfying the relations given in  page 160, \cite{KliSch-1997aa}. 
Hopf *-structure of  $U(\mathfrak{sp}(2n))$ comes from the following maps (see page $18$ and page $21$ of \cite{KliSch-1997aa}):
\[
\Delta(r)=r\otimes 1 +1 \otimes r, \quad S(r)=-r, \quad \epsilon(r)=0, \quad r=r^* \quad  \forall r \in \mathfrak{sp}(2n).
\]
Denote by $T_1$ the finite dimensional irreducible representation of $U(\mathfrak{sp}(2n))$ 
with highest weight $(1,0,\cdots ,0)$.
There exists unique nondegenerate dual pairing $\left\langle\cdot,\cdot\right\rangle$ between the
Hopf $*$-algebras $U(\mathfrak{sp}(2n))$ and $\mathcal{O}(SP(2n)/SP(2n-2))$ such that
\[
  \left\langle f,u_{l}^{k}\right\rangle = t_{kl}(f);        \hspace{1in} \mbox{for }  k= 1 \mbox{ or } 2n \mbox{ and }  1\leq l\leq 2n,\\
\]
where $t_{kl}$ is the matrix element of $T_{1}$. Using this, one can give  the algebra $\mathcal{O}(SP(2n)/SP(2n-2))$ a $U(\mathfrak{sp}(2n))$-module
structure in the following way.
\[
 f(a)=(1\otimes \langle f,.\rangle)\Delta a,
\]
where $f \in U(\mathfrak{sp}(2n))$ and $a \in \mathcal{O}(SP(2n)/SP(2n-2))$.  
We call an element $b \in \mathcal{O}(SP(2n)/SP(2n-2))$ a highest weight vector  with highest weight $(\lambda_1,\lambda_2,0,\cdots,0)$ if 
\[
 H_1(b)=(\lambda_1-\lambda_2)b, \quad H_2(b)=\lambda_2b, \quad H_i(b)= 0  \mbox{ for} \quad i\geq 2,
\]
and
\[
 E_i(b)=0 \qquad  \mbox{ for } 1\leq i \leq 2n.
\]
We will write down $\lambda_1-\lambda_2+1$ linearly independent highest weight vectors explicitly in 
terms of  $\{u_m^1, u_m^{2n}: 1 \leq m \leq 2n\} \subset \mathcal{O}(SP(2n)/SP(2n-2))$. 
Let $x= u_{2n-1}^1,y=u_{2n-1}^{2n},z=u_{2n}^1$ and $w=u_{2n}^{2n}$.   For $j \in \left\{0,1,\cdots,\lambda_1-\lambda_2\right\}$, define 
\[
 b^{(\lambda_1,\lambda_2,j)}:=z^jw^{\lambda_1-\lambda_2-j}(xw-yz)^{\lambda_2}. 
\]
\bppsn \label{new}
Let  $ \lambda_1, \lambda_2$ be two positive integers  such that
$\lambda_1 \geq \lambda_2$.
Then the set $\{b^{(\lambda_1,\lambda_2,j)}: 0 \leq j \leq \lambda_1-\lambda_2\}$ is a linearly independent set of  highest weight vectors  in the algebra
$\mathcal{O}(SP(2n)/SP(2n-2))$ with highest weight $(\lambda_1,\lambda_2,0,\cdots,0)$. 
\eppsn
\prf  It is easy to see that
\begin{displaymath}
 E_i(x)=E_i(y)=E_i(z)=E_i(w)=0 \quad \mbox{for}\quad  i >1.
\end{displaymath}
 Also, 
 \begin{displaymath}
  E_1(x)=-z, E_1(y)=-w, E_2(z)=E_2(w)=0.
 \end{displaymath}
 Further, 
 \begin{displaymath}
  H_1(x)=-x, H_1(y)=-y, H_1(z)=z, H_1(w)=w, 
 \end{displaymath}
 \begin{displaymath}
  H_2(x)=x, H_2(y)=y, H_2(z)=0, H_2(w)=0,
 \end{displaymath}
 and for $i>2$, $H_i$ maps these elements to $0$. Now using    properties of 
Hopf $*$ algebra pairing (see page $21$ of \cite{KliSch-1997aa}), one can check that  
$\{b^{(\lambda_1,\lambda_2,j)}: 0 \leq j \leq \lambda_1-\lambda_2\}$ are highest weight vectors 
with highest weight $(\lambda_1,\lambda_2,0,\cdots,0)$. 
The proof of linear independence follows from the fact that $x,y,z$ and $w$ represent projections 
or  conjugate of  projections on to different coordinates of a point in $H^{n}$.
\qed \\
Let
\[
 \Gamma=\{(\gamma_1,\gamma_2,\gamma_3): \gamma_1, \gamma_2,\gamma_3 \in \bbn, 0 \leq \gamma_2 \leq \gamma_1, 0\leq \gamma_3 \leq \gamma_1-\gamma_2\}
\]
Here first two co-ordinates  represent the highest weight and last
co-ordinate is for multiplicity. We denote by $W_{\gamma}$ the vector 
space corresponding to irreducible representation of highest weight vector $b^{\gamma}$ in the sense of Podles (see page 4, \cite{Pod-1995aa}),
by $N_{\gamma}$  the dimension of  $W_{\gamma}$
and by $\{u_{i}^{\gamma}: i\in \{1,2,\cdots N_{\gamma}\}$  a basis of $W_{\gamma}$ such that $u_{1}^{\gamma}=b^{\gamma}$. 
 Hence we can write equation (\ref{eq1}) as
\begin{IEEEeqnarray}{rCl}
\mathcal{O}(SP(2n)/SP(2n-2))= \oplus_{\gamma \in \Gamma}W_{\gamma}. \nonumber 
\end{IEEEeqnarray}
 Therefore the set $\{e_{i}^{\gamma}:=\frac{u_{i}^{\gamma}}{\|u_{i}^{\gamma}\|}: i\in \{1,2,\cdots N_{\gamma}\}, \gamma \in \Gamma\}$ 
is an orthonormal basis of $L^2(\rho)$.
Let $D$ be an equivariant Dirac operator. Then following the
arguments in propositions 5.1-5.3 leading to the statement (5.22)  in  \cite{ChaPal-2016aa}, we can assume that $D$ must be of the form 
\[
 De_{i}^{\gamma}=d^{\gamma}e_{i}^{\gamma}, \qquad i  \in \{1,2,\cdots N_{\gamma}\}, \gamma \in \Gamma.  
\]
Further assume that $(\mathcal{O}(SP(2n)/SP(2n-2)),L^2(\rho),D)$ is an equivariant spectral triple of the system $(C(SP(2n)/SP(2n-2)), SP(2n), \tau)$.
Define the set 
\[
 \Theta=\{(x,y,z,w) \in \bbr^4: 0 \leq x,y,z,w \leq 1, x^2+y^2+z^2+w^2=1\}.
\]
For $\gamma=(\gamma_1,\gamma_2,\gamma_3) \in \Lambda$, define the  function 
\[
 g^{(\gamma_1,\gamma_2,\gamma_3)} : \Theta \rightarrow \bbr
\]
sending $(x,y,z,w)$ to $z^{\gamma_3}w^{\gamma_1-\gamma_2-\gamma_3}(xw+yz)^{\gamma_2}$. Applying rotations on the co-ordinates appropriately, 
we get
\[
 \|u_{1}^{\gamma}\|=\| b^{(\gamma_1,\gamma_2,\gamma_3)}\|=\sup_{(x,y,z,w)\in \Theta } g^{(\gamma_1,\gamma_2,\gamma_3)}(x,y,z,w). 
\]
\bppsn \label{cpt}
Let $\Theta$ be a compact subset of $\bbr^n$ and $f$ and $h$ are two real valued continuous   functions define on $\Theta$.
Let $x_0 \in \Theta$ be a point such that $|f(x_0)|=\|f\|= \sup_{x \in \Theta}|f(x)| \neq 0$ and $h(x_0) \neq 0$. Then one has
\[
 \frac{\|h^mf\|}{\|h^{m+1}f\|}\leq \frac{1}{|h(x_0)|}.
\]
\eppsn
\prf  For  $m>0$, choose $x_m\in \Theta$ such that $|h^mf(x_m)|=\|h^mf\|$. Then for $m \geq 0$, we have
\[
 \|h^{m+1}f\|\geq |h^{m+1}f(x_{m})|\geq |h^mf(x_m)h(x_m)| =|h(x_m)|\|h^mf\|.
\]
Further
\[
 \|h^{m+1}f\|= |h^{m+1}f(x_{m+1})| =|h(x_{m+1})||h^mf(x_{m+1})|\leq |h(x_{m+1})|\|h^mf\|.
\]
Comparing the two inequalities, we get 
\[
 |h(x_m)|\leq |h(x_{m+1})| 
\]
Hence we have
\[
 \frac{\|h^mf\|}{\|h^{m+1}f\|}\leq\frac{1}{|h(x_m)|} \leq \frac{1}{|h(x_0)|}. 
\]
\qed

\blmma \label{estimate}
For  $(m,n) \in \bbn^2-\{0\}$, define $f_{(m,n)}: \Theta \rightarrow \bbr$ by $f_{(m,n)}(x,y,z,w)=(zw)^n(xz+yw)^m$.
Let $\theta_{(m,n)}:=(\frac{\sqrt{m}}{2\sqrt{n+m}},\frac{\sqrt{m}}{2\sqrt{n+m}},\frac{\sqrt{2n+m}}{2\sqrt{n+m}},\frac{\sqrt{2n+m}}{2\sqrt{n+m}})$.  
Then $\theta_{(m,n)} \in \Theta$ and $f_{(m,n)}(\theta_{(m,n)})=\|f_{(m,n)}\|$.
\elmma
\prf
By symmetry, we can assume without loss of generality that $x=y$ and $z=w$. Now by a straightforward calculation, one can   prove the claim.
\qed

\blmma \label{bound hwv}
Let $\epsilon_1=(1,0,0)$, $\epsilon_2=(0,1,0)$ and $\epsilon_3=(0,0,1)$. Then one has
\begin{enumerate} 
 \item
 $\sup_{\{\gamma \in \Gamma:\gamma_1=\gamma_2,\gamma_3=0\}}\frac{\|u_{1}^{\gamma}\|}{\|u_{1}^{\gamma+\epsilon_1+\epsilon_2}\|} <\infty$.
  \item
 $\sup_{\{\gamma \in \Gamma:\gamma_1-\gamma_2-2\gamma_3=0\}}\frac{\|u_{1}^{\gamma}\|}{\|u_{1}^{\gamma+2\epsilon_1+\epsilon_3}\|} <\infty$.
 \item
 $\sup_{\{\gamma \in \Gamma:\gamma_1-\gamma_2-2\gamma_3\geq 0\}}\frac{\|u_{1}^{\gamma}\|}{\|u_{1}^{\gamma+\epsilon_1}\|} <\infty$.
 \item
 $\sup_{\{\gamma \in \Gamma:\gamma_1-\gamma_2-2\gamma_3\leq 0\}}\frac{\|u_{1}^{\gamma}\|}{\|u_{1}^{\gamma+\epsilon_1+\epsilon_3}\|} <\infty$. 
\end{enumerate}
\elmma
\prf Observe that $\|g^{\gamma}\|\neq 0$ for all $\gamma \in \Gamma$.
\begin{enumerate}
 \item For $ \gamma$ with $\gamma_1=\gamma_2$ and $\gamma_3=0$, we have $g^{\gamma}=(xz+yw)^{\gamma_2}$ and 
 $g^{\gamma+\epsilon_1+\epsilon_2}=(xz+yw)^{\gamma_2+1}$ . Hence
 \begin{IEEEeqnarray}{rCl}
\sup_{\{\gamma \in \Gamma:\gamma_1=\gamma_2,\gamma_3=0\}}\frac{\|u_{1}^{\gamma}\|}{\|u_{1}^{\gamma+\epsilon_1+\epsilon_2}\|}
 &=&\sup_{\{\gamma \in \Gamma:\gamma_1=\gamma_2,\gamma_3=0\}}\frac{\|g^{\gamma}\|}{\|g^{\gamma+\epsilon_1+\epsilon_2}\|}  \nonumber \\
 &=&\frac{1}{\|g^{(1,1,0)}\| }< \infty. \nonumber 
\end{IEEEeqnarray}
\item 
For $ \gamma$ with $\gamma_1-\gamma_2-2\gamma_3=0$, we have $g^{\gamma}=(zw)^{\gamma_3}(xz+yw)^{\gamma_2}$ and 
 $g^{\gamma+\epsilon_1+\epsilon_3}=(zw)^{\gamma_3+1}(xz+yw)^{\gamma_2}$. Also, using Lemma \ref{estimate}, one can see that 
 at $(1/2,1/2,1/2,1/2)$, the function $f_{(\gamma_2,0)}(x,y,z,w)= (xz+yw)^{\gamma_2}$ takes its maximum.   Hence  by Proposition \ref{cpt}, we get
 \begin{IEEEeqnarray}{rCl}
\sup_{\{\gamma \in \Gamma:\gamma_1-\gamma_2-2\gamma_3=0\}}\frac{\|u_{1}^{\gamma}\|}{\|u_{1}^{\gamma+2\epsilon_1+\epsilon_3}\|}
 &=&\sup_{\{\gamma \in \Gamma:\gamma_1-\gamma_2-2\gamma_3=0\}}\frac{\|g^{\gamma}\|}{\|g^{\gamma+2\epsilon_1+\epsilon_3}\|} \nonumber \\
 &\leq& \frac{1}{zw(\theta_{(\gamma_2,0)})} \leq  4 < \infty. \nonumber 
\end{IEEEeqnarray}
\item 
For $ \gamma$ with $\gamma_1-\gamma_2-2\gamma_3\geq 0$, we have $g^{\gamma}=w^{\gamma_1-\gamma_2-2\gamma_3}(zw)^{\gamma_3}(xz+yw)^{\gamma_2}$.  
Hence  by  Lemma \ref{estimate} and Proposition \ref{cpt}, we get
\begin{IEEEeqnarray}{rCl}
\sup_{\{\gamma \in \Gamma:\gamma_1-\gamma_2-2\gamma_3\leq 0\}}\frac{\|u_{1}^{\gamma}\|}{\|u_{1}^{\gamma+\epsilon_1}\|}
 &=&\sup_{\{\gamma \in \Gamma:\gamma_1-\gamma_2-2\gamma_3\leq 0\}}\frac{\|g^{\gamma}\|}{\|g^{\gamma+\epsilon_1}\|}  \nonumber \\
 &\leq &\sup_{\{(\gamma_2,\gamma_3) \in \bbn^2-\{0\}\}}\frac{1}{w(\theta_{(\gamma_2,\gamma_3)}) } \leq 2 <\infty. \nonumber 
 \end{IEEEeqnarray}
\item 
For $ \gamma$ with $\gamma_1-\gamma_2-2\gamma_3\leq 0$, we have $g^{\gamma}=
z^{\gamma_2+2\gamma_3-\gamma_1}(zw)^{\gamma_1-\gamma_2-\gamma_3}(xz+yw)^{\gamma_2}$.  
Hence  by  Lemma \ref{estimate} and Proposition \ref{cpt}, we get
\begin{IEEEeqnarray}{rCl}
\sup_{\{\gamma \in \Gamma:\gamma_1-\gamma_2-2\gamma_3\geq 0\}}\frac{\|u_{1}^{\gamma}\|}{\|u_{1}^{\gamma+\epsilon_1+\epsilon_3}\|}
 &=&\sup_{\{\gamma \in \Gamma:\gamma_1-\gamma_2-2\gamma_3\geq 0\}}\frac{\|g^{\gamma}\|}{\|g^{\gamma+\epsilon_1+\epsilon_3}\|}  \nonumber \\
 &\leq &\sup_{\{(\gamma_2,\gamma_3) \in \bbn^2-\{0\}\}}\frac{1}{z(\theta_{(\gamma_2,\gamma_1-\gamma_2-\gamma_3)}) } \leq 2 <\infty. \nonumber
\end{IEEEeqnarray}
\end{enumerate}
\qed \\
In the following lemma, we establish some estimate on the growth of $d^{\gamma}$.
\blmma \label{bound}
Let $\epsilon_1=(1,0,0)$, $\epsilon_2=(0,1,0)$ and $\epsilon_3=(0,0,1)$. Then one has
\begin{enumerate} 
 \item
 $\sup_{\{\gamma \in \Gamma:\gamma_1=\gamma_2,\gamma_3=0\}}|d^{\gamma+\epsilon_1+\epsilon_2}-d^{\gamma}| <\infty$.
  \item
 $\sup_{\{\gamma \in \Gamma:\gamma_1-\gamma_2-2\gamma_3=0\}}|d^{\gamma+2\epsilon_1+\epsilon_3}-d^{\gamma}| <\infty$.
 \item
 $\sup_{\{\gamma \in \Gamma:\gamma_1-\gamma_2-2\gamma_3\geq 0\}}|d^{\gamma+\epsilon_1}-d^{\gamma}| <\infty$.
 \item
 $\sup_{\{\gamma \in \Gamma:\gamma_1-\gamma_2-2\gamma_3\leq 0\}}|d^{\gamma+\epsilon_1+\epsilon_3}-d^{\gamma}|<\infty$.
\end{enumerate}
\elmma
\prf
Let $\gamma \in \Gamma$ such that $\gamma_1=\gamma_2$ and $\gamma_3=0$. Note that 
\begin{IEEEeqnarray}{rCl}
\|[D, (xw-yz)]u_{1}^{\gamma}\|&=&\|[D, (xw-yz)](xw-yz)^{\gamma_2}\| \nonumber  \\
&=&\|(d^{\gamma}-d^{\gamma+\epsilon_1+\epsilon_2})u_{1}^{\gamma+\epsilon_1+\epsilon_2}\|. \nonumber 
\end{IEEEeqnarray}
Since $[D, (xw-yz)]$ is a bounded operator,  we have 
\begin{IEEEeqnarray}{rCl}
\sup_{\{\gamma \in \Gamma:\gamma_1=\gamma_2,\gamma_3=0\}}|d^{\gamma}-d^{\gamma+\epsilon_1+\epsilon_2}|&=&
\sup_{\{\gamma \in \Gamma:\gamma_1=\gamma_2,\gamma_3=0\}}\frac{\|[D, (xw-yz)]u_{1}^{\gamma}\|}{\|u_{1}^{\gamma+\epsilon_1+\epsilon_2}\|} \nonumber \\
&\leq & \sup_{\{\gamma \in \Gamma:\gamma_1=\gamma_2,\gamma_3=0\}} \|[D, (xw-yz)]\| \frac{\|u_{1}^{\gamma}\|}{\|u_{1}^{\gamma+\epsilon_1+\epsilon_2}\|} \nonumber \\
& < & \infty \hspace{0.3in} ( \mbox{ by part (1) of the Lemma } \ref{bound hwv})\nonumber
\end{IEEEeqnarray}
Using Lemma \ref{bound hwv} and the  fact that  $[D, zw]$, $[D, z]$ and $[D, w]$ are bounded operators, other parts of the claim follow similarly.
\qed \\
Let $c>0$ be an upper bound in all the four inequalities  of the Lemma \ref{bound}.  Let $\mathcal{G}$ be a graph with vertex set $\Gamma$ and edge set
$\{(\gamma, \gamma^{'}): |d^{\gamma}-d^{\gamma^{'}}|< c\}$. The following lemma says that $\mathcal{G}$ is a connected graph.
\blmma \label{path}
Let $\gamma \in \Gamma$. Then there is a path in $\mathcal{G}$ joining $(0,0,0)$ and  $\gamma$ and of length  less than or equal to $\gamma_1$.
\elmma
\prf 
If $\gamma_1-\gamma_2-2\gamma_3\geq 0$, then one possible path would be as follows.
\begin{IEEEeqnarray}{lCl}
(0,0,0)\rightarrow (1,1,0)\rightarrow (2,2,0)\rightarrow \cdots \rightarrow (\gamma_2,\gamma_2,0) \nonumber \\
\hspace{2in} \mbox{(by part(1) of the Lemma } \ref{bound}) \nonumber \\
 (\gamma_2,\gamma_2,0) \rightarrow (\gamma_2+2,\gamma_2,1) \rightarrow \cdots \rightarrow
 (\gamma_2+2\gamma_3,\gamma_2,\gamma_3)  \nonumber \\
  \hspace{2in} \mbox{ (by part(2) of the  Lemma } \ref{bound}) \nonumber \\
 (\gamma_2+2\gamma_3,\gamma_2,\gamma_3)  \rightarrow 
 (\gamma_2+2\gamma_3+1,\gamma_2,\gamma_3)  \rightarrow \cdots (\gamma_1,\gamma_2,\gamma_3) \nonumber \\
\hspace{2in} \mbox{ (by part(3) of the  Lemma } \ref{bound}).   \nonumber
\end{IEEEeqnarray}
If $\gamma_1-\gamma_2-2\gamma_3\leq 0$, then one possible path would be as follows.
\begin{IEEEeqnarray}{lCl}
(0,0,0)\rightarrow (1,1,0)\rightarrow (2,2,0)\rightarrow \cdots \rightarrow (\gamma_2,\gamma_2,0) \nonumber \\
\hspace{2in} \mbox{(by part(1) of the Lemma } \ref{bound}) \nonumber \\
 (\gamma_2,\gamma_2,0) \rightarrow (\gamma_2+2,\gamma_2,1) \rightarrow \cdots \rightarrow
 (2\gamma_1-\gamma_2-2\gamma_3,\gamma_2,\gamma_1-\gamma_2-\gamma_3)  \nonumber \\
  \hspace{2in} \mbox{ (by part(2) of the  Lemma } \ref{bound}) \nonumber \\
 (2\gamma_1-\gamma_2-2\gamma_3,\gamma_2,\gamma_1-\gamma_2-\gamma_3) \rightarrow 
 (2\gamma_1-\gamma_2-2\gamma_3+1,\gamma_2,\gamma_1-\gamma_2-\gamma_3+1)  \nonumber \\  
 \rightarrow \cdots  \rightarrow  (\gamma_1,\gamma_2,\gamma_3) 
\hspace{.7in} \mbox{ (by part(4) of the  Lemma } \ref{bound}).   \nonumber
\end{IEEEeqnarray}
Moreover, the length of the paths in both cases are less than $\gamma_1$ as in 
each step the increament in the first coordinate is at least one. This settles the claim.
\qed 
\blmma \label{bounded growth}
Let $D: e_{i}^{\gamma} \mapsto d^{\gamma}e_{i}^{\gamma}$ be an operator acting on the Hilbert space $L^2(\rho)$ such that the triple
$(L^2(\rho),\pi_{\rho},D)$ is an equivariant spectral triple  of the system $(C(SP(2n)/SP(2n-2)), SP(2n), \tau)$. Then we have 
\[
 d^{\gamma}=O(\gamma_1).
\]
\elmma
\prf
It follows from  Lemma \ref{path}. 
\qed 
\blmma \label{bounded leap}
For $1\leq m \leq 2n$ and $l=1$ or $2n$, one has 
\[
 u_m^lu_{i}^{\gamma} \subset \mbox{span}\{ u_{i}^{\beta}: \gamma_1-1 \leq \beta_1 \leq \gamma_1+1 \}
\]
\elmma
\prf
Let $e_i=(0,0,\cdots, 0, \underbrace{1}_{i^{th}-\mbox{place}},0,\cdots 0)$. Then from equation ((14), page $210$, \cite{KliSch-1997aa}), we get 
\[
 u_{(\gamma_1,\gamma_2,0,0,\cdots ,0)} \otimes u_{(1,0,\cdots ,0)} = \oplus_{i=1}^n u_{(\gamma_1,\gamma_2,0,0,\cdots ,0) + e_i} \oplus \oplus_{i=1}^n u_{(\gamma_1,\gamma_2,0,0,\cdots ,0) - e_i}
\]
Hence $u_m^lu_{i}^{\gamma}$ is in the span of matrix entries of the irreducible representations of  highest weight $(\beta_1, \beta_2, \cdots, \beta_n)$ such that $\beta_1= \gamma_1$ or $\gamma_1 \pm 1$. 
Since  $u_m^lu_{i}^{\gamma} \in \mathcal{O}(SP(2n)/SP(2n-2))$ and   $\{u_{i}^{\gamma}: \gamma \in \Gamma\}$ is a basis of $\mathcal{O}(SP(2n)/SP(2n-2))$,  we get the claim.
\qed 
\bthm \label{optimal}
Let $D_{eq}$ be the Dirac operator  $e_{i}^{\gamma} \mapsto \gamma_1e_{i}^{\gamma}$ acting on the Hilbert space $L^2(\rho)$. Then  the triple
$(\mathcal{O}(SP(2n)/SP(2n-2)),L^2(\rho),D)$ is a $(4n-1)$-summable equivariant spectral triple  of the  system $(C(SP(2n)/SP(2n-2)), SP(2n),\tau)$.
The operator $D_{eq}$ is optimal, i.e. if $D$ is any equivariant Dirac operator of the $C^*$-dynamical system $(C(SP(2n)/SP(2n-2)), SP(2n), \tau)$ acting on $L^2(\rho)$ then there exist positive reals $a$ and $b$ such that
\[
 |D|\leq a|D_{eq}|+b.
\]
\ethm 
\prf Clearly $D_{eq}$ is a selfadjoint operator with compact resolvent.
That $D_{eq}$ has bounded commutators with the generators $\left\{u_m^1,u_m^{2n} : m \in \{1,2,\cdots 2n\}\right\}$ of $\mathcal{O}(SP(2n)/SP(2n-2))$
follows from Lemma \ref{bounded leap}. This proves that the triple $(\mathcal{O}(SP(2n)/SP(2n-2)), L^2(\rho), D)$ is an 
equivariant spectral triple  of the  system $(C(SP(2n)/SP(2n-2)), SP(2n), \tau)$.
From Weyl dimension formula, we have
\[
 N_{\gamma}=O(\gamma_1^{2n-1}\gamma_2^{2n-3}).
\]
This along with the fact that $0 \leq \gamma_2 \leq \gamma_1$ and  $0\leq \gamma_3 \leq \gamma_1-\gamma_2$ shows that $D$ is $(4n-1)$-summable. Optimality follows from Lemma \ref{bounded growth}.
\qed 
\bthm \label{spectral}
Spectral dimension of the  quaternion spheres $SP(2n)/SP(2n-2)$ is $4n-1$.
\ethm
\prf It is a direct consequence of Theorem \ref{optimal}.
\qed

\noindent{\sc Bipul Saurabh} (\texttt{saurabhbipul2@gmail.com})\\
         {\footnotesize Harish Chandra Research Institute , Chhatnag Road, Jhunsi,
Allahabad 211019,  INDIA}

\end{document}